\documentclass[11pt]{amsart}
\usepackage{amsmath,amssymb,enumerate}
\usepackage{latexsym}
\usepackage{amsfonts}
\usepackage {amsbsy}
\usepackage {amsmath}
\usepackage{amssymb}
\usepackage{enumerate}
\usepackage{ bbm}

\newtheorem{theorem}{Theorem}
\newtheorem{lemma}[theorem]{Lemma}
\newtheorem{corollary}[theorem]{Corollary}

\theoremstyle{definition}
\newtheorem{definition}[theorem]{Definition}

\numberwithin{theorem}{section}
\numberwithin{equation}{section}

\theoremstyle{definition}

 \numberwithin{equation}{section} 
 \numberwithin{figure}{section} 
 \theoremstyle{plain}    
 \theoremstyle{plain}    
 \theoremstyle{remark}
 \theoremstyle{remark}
 \theoremstyle{definition}
\theoremstyle{plain}  
\theoremstyle{plain}

\theoremstyle{definition}

\numberwithin{equation}{section}

\newcommand{\N}{\mathbb{N}}
\newcommand{\R}{\mathbb{R}}
\newcommand{\C}{\mathbb{C}}

\makeatletter
\@namedef{subjclassname@2010}{%
  \textup{2010} Mathematics Subject Classification}
\makeatother

\begin{document}

\keywords{Complex Monge-Ampère operator, Dirichlet Eigenvalue, Iterative method.}

\subjclass[2010]{31C45, 32U15, 32U40, 32W20, 35J66, 35J96}

\title[An Iterative Approach ]{An Iterative Approach to the Complex Monge-Amp\`ere Eigenvalue Problem}

\author{Ahmed Zeriahi}

\address{Institut de Mathématiques de Toulouse; UMR 5219, Université de Toulouse; CNRS, UPS, 118 route de Narbonne, F-31062 Toulouse Cedex 9, France}

\email{ahmed.zeriahi@math.univ-toulouse.fr}

\date{\today}

\maketitle

\begin{center}
{\it In honor to Wies\l aw Plesniak \\ on his $80^{\text{th}}$ birthday}
\end{center}

\begin{abstract} We present an iterative approach to approximate the solution to the Dirichlet  complex Monge-Ampère eigenvalue problem on a bounded strictly pseudoconvex domain in $\C^n$. This approach is inspired by a similar approach initiated by F. Abedin, J. Kitagawa who considered the real Monge-Ampère operator on a strictly convex domain in $\R^N$.

\smallskip

This work is based on recent results obtained by P. Badiane and the author on the existence and uniqueness of the solution to the Dirichlet eigenvalue problem for the complex Monge-Ampère operator.

\smallskip

However, the iterative approach does not require the a priori knowledge of the first eigenvalue but it provides an  {effective scheme} for approximating it, as well as the associated eigenfunction.
\end{abstract}

\section{Introduction}

\subsection{ The complex Monge-Amp\`ere eigenvalue problem }

Let $\Omega \Subset  \C^n $ be a bounded strictly pseudoconvex domain.
Let $\mu$ be a smooth non degenerate volume form (measure) on $\bar \Omega$ i.e. $ \mu  := f dV$, where  $0 < f  \in C^{\infty}(\bar \Omega)$ and $dV$ is  the Euclidean volume form on $\C^n$.

\smallskip
\smallskip

{\bf The eigenvalue problem  :}   Find a pair $(\lambda,u)$, where  $\lambda > 0$ is a constant and  $u \in  {PSH} (\Omega)\cap C^{0}(\bar \Omega),$ satisfying the following properties:
\begin{equation}\label{eq:VP-MA}
\left\{\begin{array}{lcll} 
 (dd^c u)^n = (- \lambda \, u)^n f dV  &\hbox{on}&  \Omega,\\
  u = 0 &  \hbox{in}&  \partial \Omega, \\
 u \not \equiv 0 & \hbox{in}& \Omega.&
\end{array}\right.
\end{equation} 
Here we use the standard differential operators $d = \partial + \bar \partial$, $d^c := (i \slash 2) (\bar \partial - \partial)$ so that $dd^c = i \partial \bar \partial$ and the volume form $d V$ is defined by
$d V = \beta^n$, where   $\beta:= dd^c\vert z \vert^2$ is the standard Kähler form on $\C^n$.

Recall that $(dd^c u)^n$ is the Monge-Ampère measure (current) associated to the function  $u \in  {PSH} (\Omega)\cap C^{0}(\bar \Omega),$ as defined by E. Bedford and B.A. Taylor (see Section 2). Observe that when $u \in PSH(\Omega) \cap C^{2}(\Omega)$, its Monge-Ampère measure is given by the following formula : 
$$
(dd^c u)^n = \text{det} \, \left( \frac{\partial^2 u}{\partial z_j \partial \bar z_k}\right) d V,
$$
considered as a Radon measure on $\Omega$ with density with respect the  Lebesgue measure defined by the  Euclidean volume form $d V$ on $C^n$.

This problem has been solved in \cite{BZ23}. Let us recall the main results of this paper for later reference.
 
 \smallskip
 \smallskip
 
{\bf Theorem A} \label{thm:Variational1} {\it Under the assumptions above, the following properties hold:

\smallskip

\noindent $(i)$  there exists a number $\lambda_1 > 0$ and a function $u_1 \in PSH (\Omega) \cap  C^{1,\bar 1}(\bar \Omega) \cap C^{\infty}(\Omega)$ such that  $(\lambda_1,u_1)$ is a solution to the eigenvalue problem \eqref{eq:VP-MA}. 

\smallskip

\noindent $(ii)$ if $\tilde{\lambda}_1 > 0$,  $\tilde{u}_1 \in PSH (\Omega) \cap C^{\alpha}(\bar \Omega)$   (with $\alpha \in ]0,1[$) and $(\tilde{u}_1,\tilde{\lambda}_1)$  is another (weak) solution to the eigenvalue problem \eqref{eq:VP-MA}, 
 then $\tilde{\lambda}_1 = \lambda_1$ and there exists a constant $a >0$ such that  $\tilde{u}_1 = a \, u_1$. }

\smallskip

Here $u \in C^{1,\bar 1}(\bar \Omega)$ means that $u \in C^{0,1} (\bar \Omega)$ and $\Delta u \in L^{\infty}(\Omega)$. In particular $u \in C^{1,\alpha}(\bar \Omega)$ for any $0 < \alpha < 1$.

\smallskip

The number $\lambda_1 = \lambda_1(\Omega,\mu)$ is called the {\it first Dirichlet eigenvalue} of the Complex Monge-Amp\`ere Operator associated to the measure space $(\Omega,\mu)$, where  $\mu = f d V$.

\smallskip
\smallskip

The second main result of \cite{BZ23} is as follows.

\smallskip

{\bf Theorem B} {\it The first eigenvalue $\lambda_1 = \lambda_1(\Omega,\mu)$ satisfies the following Rayleigh quotient type formula : 
  \begin{equation} \label{eqR1}
  \lambda_1^n = \frac{E(\varphi_1)}{I_\mu (\varphi_1)} = \inf \left\{\frac{E(\phi)}{I_\mu (\phi)} ; \phi \in \mathcal E^1 (\Omega), \phi \not \equiv 0 \right\}\cdot
 \end{equation}  }

\smallskip

\smallskip

Recall that $\mathcal E^1 (\Omega)$ is the Cegrell class of plurisubharmonic   functions in $\Omega$ with finite Monge-Amp\`ere energy (see Section 2).

\smallskip

The energy functional 
$E : \mathcal E^1 (\Omega) \longrightarrow \R^+$ is defined  by the formula :  
$$
  \mathcal E^1 (\Omega) \ni \phi \longmapsto E(\phi) :=  \frac{1}{ n+1} \int_\Omega (-\phi) (dd^c \phi)^n,  \, \, \, \, 
$$
where $(dd^c \phi)^n$ is the Monge-Ampère measure associated to $\phi$ (see Section 2).
\smallskip

The integral functional  $I_\mu : \mathcal E^1 (\Omega) \longrightarrow \R^+$ is defined by the formula :  
$$
 \mathcal E^1 (\Omega) \ni \phi \longmapsto I_\mu(\phi) :=  \frac{1}{ n+1}  \int_\Omega (-\phi)^{n +1} d \mu,
$$
provided that $\mathcal E^1 (\Omega) \subset L^{n+1}(\Omega,\mu)$ as we will show.

\smallskip
\smallskip

\noindent{\bf Previous known results :}
\begin{itemize}
\item  It's well known that  a strictly elliptic symmetric second order operator (e.g. the Laplace operator $- \Delta_{\R^n}$) acting on smooth functions in a bounded domain with smooth boundary in $\R^n$,  admits infinitely many Dirichlet eigenvalues $0 < \lambda_1 < \lambda_2 < \cdots < \lambda_k \to + \infty$. 

 Moreover, among all these eigenvalues, the first one $\lambda_1$ is the only one for which there exists an eigenfunction with a constant sign.

\smallskip

\item For the real Monge-Amp\`ere operator  acting on convex functions on a bounded strictly convex domain in $\R^n$, the existence and uniqueness of the first Dirichlet  eigenvalue is due to P.L. Lions \cite{Lions86}. 

\smallskip

\item The Rayleigh quotient formula was proved by K. Tso using parabolic methods \cite{Tso90}.
\end{itemize}

\smallskip

\noindent{\bf Comments :}  1) We provided in \cite{BZ23} two different approaches: 
\begin{itemize}
\item {The PDE approach} following P.-L. Lions \cite{Lions86} original approach in the real case with specific a priori estimates for solutions to some degenerate complex Monge-Amp\`ere equations.
 
\smallskip

\item {The Variational approach}  in the spirit of the variational approach to the resolution of complex Monge-Amp\`ere equations introduced in our earlier work \cite{BBGZ13} (see also \cite{ACC12}).
 
\smallskip

\noindent 2) Let us also mention the following fact.
 Recently using our strategy,  J. Chu, Y.Liu, N. McCleerey \cite{CLMcC24} have extended our results to the Hessian operators and showed that the eigenfunctions are actually in $C^{1, 1}(\bar \Omega)$. 
\end{itemize}

\smallskip

Let us describe briefly these two approaches. 
 
\smallskip 

\subsection{ The PDE approach.} This method is inspired by the work of Lions  \cite{Lions86}. We first use an observation due to B. Gaveau \cite{Gav77}. 

Let $\mathcal A(\Omega)$ be the set of hermitian $n \times n$ matrices $a =(a_{j \bar k})$
with smooth entries on $\Omega$ such that $ \text{det} \, a \geq 1$ in $\Omega$. 

 To each $a \in \mathcal A (\Omega)$ we associate the linear second order elliptic  operator defined on smooth functions $u \in C^2(\Omega)$ by the following formula:
$$
L_a u := \frac{1}{n} \sum_{1\leq j, k \leq n} a_{j \bar k}  u_{i \bar k}, \, \, \text{where} \, \, \, u_{j \bar k} := \frac{\partial^2 u}{\partial z_j \partial \bar z_k}\cdot
$$

Then for a function $u \in PSH(\Omega) \cap C^2(\Omega)$, we have the following formula:
$$
\left[\text{det} \, (u_{j \bar k})\right]^{1\slash n} = \inf \{L_a u ; a \in \mathcal A(\Omega)\}\cdot
$$

For each $a \in \mathcal A(\Omega)$, we denote by $\gamma_1(a) > 0$ the first eigenvalue of the linear operator $-L_a$ and set
$$
\lambda_1 := \inf \{\gamma_1(a) \, ; \, a \in \mathcal A(\Omega)\}.
$$

On the other hand, consider the following Dirichlet problem for the complex Monge-Amp\`ere equation :
\begin{equation}\label{eq:Approximate}
\left\{\begin{array}{lcll} 
 (dd^c u)^n = (1 - \lambda \, u)^n f dV  &\hbox{on}&  \Omega,\\
  u = 0 &  \hbox{in}&  \partial \Omega,
\end{array}\right.
\end{equation} 
where $\lambda > 0$ is a parameter.

Note that as the RHS is smooth, non degenerate but  decreasing in $u$, the existence result of  L. Caffarelli, J.J. Kohn, L. Nirenberg, J. Spruck \cite{CKNS85} does not apply. Actually neither the existence nor the uniqueness of a solution to the equation \eqref{eq:Approximate} is guaranteed.

\smallskip

Let $\Lambda \subset \R^+$ be the set of real numbers $\lambda \geq 0$ such that the problem \eqref{eq:Approximate} has a smooth solution. 
By the fundamental  result of  \cite{CKNS85}, there exists a unique   $u_0 \in PSH(\Omega) \cap C^{\infty}(\bar \Omega)$ such that $(dd^c u_0)^n = f d V$ and ${u_0}_{\mid \partial \Omega} \equiv 0$. Hence  $0 \in \Lambda$, and  $\Lambda \neq \emptyset$.

\smallskip
\smallskip

\noindent{\bf Claim 1 : } {$\sup \Lambda = \lambda_1$. 

\smallskip

To prove this claim we need to establish a new existence theorem for the Dirichlet problem. 
If the problem \eqref{eq:Approximate} admits a smooth subsolution $\underline u$ and a weak supersolution $\bar u$  such that $\underline u \leq \bar u$ in $\Omega$, then it admits a smooth solution $u$ such that $ \underline u \leq u \leq \bar u$. 

 To this end, we had to establish new a priori estimates for solutions to degenerate complex  Monge-Amp\`ere equations \eqref{eq:Approximate} (see \cite{BZ23}).

\subsection{The Variational approach.} We consider the minimization problem for the Rayleigh quotient :
$$
\eta_1^n := \inf \left\{\frac{E(\psi)}{I_\mu(\phi)} \, ; \, \phi \in \mathcal E^1(\Omega) \setminus \{0\}\right\}.
$$

 Then we show that any normalized minimizing sequence 
$(\psi_j)_{j \in \N}$ in $\mathcal E^1(\Omega) \setminus \{0\}$ of the Rayleigh quotient admits a subsequence which converges in $L^1_{loc}(\Omega)$ to a function $\psi \in \mathcal E^1(\Omega) \setminus \{0\}$  which satisfies
$E(\psi) = \eta_1^n I_\mu(\psi).$

Therefore $\psi \not \equiv 0$ is a minimizer of the functional defined on $\mathcal E^1(\Omega)$ by
$$
 F_\mu (\phi) := E(\phi) - \eta_1^n I_\mu (\phi).
$$

To prove that $(\eta_1,\psi)$ is a solution to the eigenvalue problem, we observe that
the Euler-Lagrange equation of the  functional $F_\mu$ is formally given by the equation
$ F_\mu'(\psi) = -  (dd^c \psi)^n +  \eta_1^n (-\psi)^n \mu = 0$. 
Using an envelope trick inspired by \cite{BBGZ13}, we show indeed that $\psi$ is a solution to the complex-Monge-Ampère equation 
$$
(dd^c \psi)^n =   \eta_1^n (-\psi)^n \mu.
$$
From Kolodziej type a priori estimates we conclude that $\psi \in C^\alpha(\bar \Omega)$ for any $\alpha \in ]0,1 \slash (2nq+1)[$  (see next subsection 3.4).
Uniqueness is done in two steps using an argument similar to that of \cite{Lions86} based on a reduction to the linear case and an application of Hopf's lemma.

\section{Statement of the Main Theorem}
  
  \subsection{The Rayleigh quotient} To define the Rayleigh quotient we  need the following properties which will be proved in the next section (see Lemma \ref{lem:integrability}).
  \smallskip
  
  \noindent{\bf Claim 2 :} {\it Under the assumptions of the previous section, we have the following properties:
  \begin{itemize}
  \item $\mathcal E^1 (\Omega) \subset L^{n+1}(\Omega,\mu)   \, \text{and} \, \,$
  \item $   \phi \in \mathcal E^1 (\Omega) \setminus \{0\}  \Longrightarrow 0 < I_\mu(\phi) < +\infty.$
  \end{itemize} }  
  \smallskip
  
   Thanks to this result, we can define the {\it Rayleigh quotient} of any $\phi \in \mathcal E^1 (\Omega)\setminus \{0\}$ : 
  $$
  R (\phi) = R_\mu (\phi) =  := \frac{E(\phi)}{I_\mu (\phi)} = \frac{\int_\Omega(-\phi) (dd^c \phi)^n}{\int_\Omega (-\phi)^{n+1} \mu}\cdot
  $$
  
  \noindent{\it Observation  : }  if $(\lambda,u)$ is a non trivial solution to the equation \eqref{eq:VP-MA} i.e.
  $$
  (dd^c u)^n =(-\lambda u)^n \mu,
  $$
  
   then
  $\lambda^n = R(u)$, hence $u$ is a {\it non-trivial  solution} to the Dirichlet problem 
  \begin{equation}\label{eq:VP-MA1}
\left\{\begin{array}{lcll} 
 (dd^c u)^n = R(u)  (- u)^n \mu &\hbox{on}&  \Omega,\\
  u = 0 &  \hbox{in}&  \partial \Omega.
\end{array}\right.
\end{equation}
 
 This is a fixed point problem for the inverse of the Complex Monge-Amp\`ere operator.
 
  More precisely,  we denote by $PSH_0(\Omega)$ the set of bounded plurisubharmonic functions in $\Omega$ with boundary values $0$ and let $\mathcal P^*_0(\Omega) := PSH_0(\Omega) \cap C^0(\bar \Omega) \setminus \{0\}$.   
  
  Let $ \phi \in  \mathcal P^*_0(\Omega)$. Since the function $R(\phi)  (- \phi)^n f $ is continuous in $\bar \Omega$, it follows from  a result of Bedford and Taylor \cite{BT76}, that there exists  a  unique solution $\psi := T(\phi) \in PSH_0(\Omega) \cap C^0(\bar \Omega) $ to the Dirichlet problem
 \begin{equation}\label{eq:MA}
\left\{\begin{array}{lcll} 
 (dd^c \psi)^n = R(\phi)  (- \phi)^n f d V &\hbox{on}&  \Omega,\\
  \psi = 0 &  \hbox{in}&  \partial \Omega.
\end{array}\right.
\end{equation}
Then it is clear that $ \psi \not \equiv 0$, hence $\psi \in \mathcal P^*_0(\Omega)$. This defines an operator $T : \mathcal P^*_0(\Omega) \longrightarrow \mathcal P^*_0(\Omega)$.

\smallskip

\noindent{ \it Reformulation of the problem :} A function $u \in \mathcal P^*_0(\Omega)$ is an eigenfunction for the eigenvalue problem, if and only if $u$ is a {\it non-trivial fixed point} of the  operator $T : \mathcal P^*_0(\Omega) \longrightarrow \mathcal P^*_0(\Omega)$.

\smallskip

We already know by Theorem A, that such eigenfunction exists, but we want to find an {\it efficient scheme} to approximate it.

\smallskip

\smallskip

\subsection{Construction of the  iterative sequence } 
 The goal is to find an initial function $u_0 \in \mathcal P_0^*(\Omega)$ so that  the sequence of its iterates $k \longmapsto u_k := T^k (u_0)$ converges  uniformly in $\bar \Omega$ to a function $ u$, and  prove that $u \in \mathcal P_0^*(\Omega)$ is a (non trivial) fixed point of $T$.
  
   \smallskip
   
  Indeed, let $u_0 \in PSH(\Omega) \cap C^{0} (\bar \Omega)$ be such that :
  \begin{equation}\label{eq:SubSol}
\left\{\begin{array}{lcll} 
 (dd^c u_0)^n \geq  f dV, \, \, &\hbox{on}&  \Omega,&\\
  {u_0}_{\mid \partial \Omega} \leq 0. 
\end{array}\right.
\end{equation} 
  For example the function $u_0 := A \rho$ with $A >> 1$ satisfies the previous requirement. Here $\rho$ is  strictly plurisubharmonic defining function of $\Omega$.
  
   Observe that  by Lemma \ref{lem:integrability}, we have $ 0 < R(u_0) < +\infty$. 
  Then starting from $u_0$ we  define by induction a  sequence $(u_k)_{k \in \N}$ such for each $k \in \N$, $u_k \in PSH(\Omega) \cap C^{0}(\bar \Omega)$, $0< R(u_k) < \infty$ and $u_{k+1} = T (u_k) \in PSH(\Omega) \cap C^0 (\bar \Omega)$,  i.e. $ u_{k+1}$ is the unique  solution to the following Dirichlet problem:
 \begin{equation}\label{eq:DP}
\left\{\begin{array}{lcll} 
 (dd^c u_{k+1})^n = R(u_k) (-u_k)^n f dV  &\hbox{on}&  \Omega,&\\
  u_{k +1} = 0 &  \hbox{in}&  \partial \Omega.&  
\end{array}\right.
\end{equation} 
Indeed, since for a fixed $k \in \N$, the right hand side $ f_k := R(u_k) (-u_k)^n f  \in C^{0} (\bar \Omega)$ and $f_k > 0$ in $\Omega$,  by a  theorem of Bedford-Taylor \cite{BT76}, the Dirichlet Problem \eqref{eq:DP} admits a unique solution $u_{k+1} \in \mathcal P_0^*(\Omega)$.

 \smallskip

Our main result is the following.

\smallskip
\smallskip

 \noindent{\bf Main Theorem} {\it  Let $u_0 \in  PSH(\Omega) \cap C^{0,1} (\bar \Omega)$ be a function satisfying \eqref{eq:SubSol}. 
   Then the  sequence $(u_k)_{k \in \N}$ converges uniformly in $\bar \Omega$ to a function $ u \in \mathcal P_0^*(\Omega)$, satisfying the following properties:
  
  $(1)$ $\lim_{k \to + \infty} R(u_k) = \lambda_1^n = R(u)$, where $\lambda_1= \lambda_1(\Omega,\mu)$
  is  the  first eigenvalue of the complex Monge-Amp\`ere operator  associated to $(\Omega,\mu)$.
  
  $(2)$ $ u \in PSH(\Omega) \cap C^{\alpha}(\bar \Omega)$ for some $\alpha \in ]0,1[$ and is an eigenfunction  associated  to the eigenvalue  $\lambda_1$ i.e. $(dd^c u)^n = (- \lambda_1 u)^n f d V $ on $  \Omega$ and $u_{\mid \partial \Omega} \equiv 0$.}
  
  \smallskip
  
  This result is inspired by a similar result due to F. Abedin and J. Kitagawa \cite{AK20} who considered  the real Monge-Amp\`ere operator in a bounded  convex domain in $\R^n$.
  However the method of proof we use here is different and relies on deep results from Pluripotential Theory as we will see.

\smallskip

Observe that this kind of approximation scheme is related to the so called  {\it inverse iteration methods} used in Linear Algebra  to approximate the smallest and the largest eigenvalues of linear symmetric operators.
 
 \section{Preliminary results}

\subsection{The Monge-Amp\`ere energy functional}
We recall the following well known  definitions.
  \begin{definition} \label{def:hyperconvex}
1. We say that a bounded domain $\Omega \Subset \C^n$ is hyperconvex if there exists  a bounded continuous plurisubharmonic  function $\rho: \Omega \to ]-\infty,0[$ which is exhaustive i.e. for any $c < 0$, $\{z \in \Omega \, ; \, \rho(z) < c\}\Subset \Omega$. 

\smallskip

2. We say that a domain $\Omega \Subset \C^n$ is strictly pseudoconvex if $\Omega$ admits a smooth defining function $\rho$ which is strictly plurisubharmonic in a neighborhood of $\bar \Omega$ and satisfies $\vert \nabla \rho \vert >0$ pointwise in $\partial \Omega=\{\rho=0\}$. 
In this case we can choose $\rho$ so that
\begin{equation} \label{eq:rho}
dd^c\rho \geq \beta,
\end{equation}
pointwise in $\Omega$. 
\end{definition}

In this section we assume that $\Omega \Subset \C^n$ is hyperconvex. Let us define some convex classes of singular plurisubharmonic functions in $\Omega$ suitable for the variational approach.  

Recall first that by E. Bedford and B.A. Taylor \cite{BT76,BT82}, the complex Monge-Ampère operator $(dd^c \cdot)^n$ is well defined on the class of locally bounded plurisubharmonic functions in $\Omega$. For such a function $u \in PSH(\Omega) \cap L_{loc}^{\infty} (\Omega)$, $(dd^c u)^n$ is a positive current of bidegree $(n,n)$ on $\Omega$, which will be identified to a Radon measure on $\Omega$. This Borel measure can be obtained as the unique limit of the Monge-Amp\`ere measures $\{(dd^c u_j)^n\}_{j \in \N}$ of any decreasing sequence of smooth plurisubharmonic functions $(u_j)_{j \in \N}$ converging to $u$ in $\Omega$. 

 By the classical representation theorem of F. Riesz, this Radon measure extends (uniquely) to a positive Borel measure on $\Omega$ (with locally finite mass), called the complex Monge-Ampère measure of $u$ and still denoted by $(dd^c u)^n$. 

 Following Urban Cegrell \cite{Ceg98}, we define the class $\mathcal{E}^0 (\Omega)$ of plurisubharmonic test functions on $\Omega$. This is the set of bounded plurisubharmonic functions $\phi$ in $\Omega$ with boundary values $0$ such that $\int_\Omega (dd^c \phi)^n < + \infty$. Then we define 
$\mathcal{E}^1 (\Omega)$ as the set of plurisubharmonic functions $u$ in $\Omega$ such that there exists a decreasing sequence $(u_j)_{j \in \N}$ in the class $\mathcal{E}^0 (\Omega)$  satisfying $u = \lim_j u_j $ in $\Omega$ and $\sup_j \int_\Omega (-u_j) (dd^c u_j)^n < + \infty$. 
It is easy to see from the definition that  $\mathcal{E}^1 (\Omega)$ is a convex cone in $L^1_{loc}(\Omega)$.  

It turns out that the complex Monge-Amp\`ere operator extends to the class $\mathcal{E}^1 (\Omega)$ and is continuous under monotone limits in $\mathcal{E}^1 (\Omega)$. Moreover if $u \in \mathcal{E}^1 (\Omega)$, then $\int_\Omega (-u) (dd^c u)^n < + \infty$ (see  \cite{Ceg98}).

The Monge-Amp\`ere energy functional is defined on the space $\mathcal{E}^1 (\Omega)$  as follows:  for $\phi \in \mathcal{E}^1 (\Omega)$,
\begin{equation} \label{eq:energyF}
E (\phi) :=  \frac{1}{n+1} \int_\Omega (-\phi) (dd^c \phi)^n.
\end{equation}
 It was proved in  \cite{BBGZ13} in the compact K\"ahler setting that the functional $E$ is lower semi-continuous and  $-E$ is a primitive of the complex Monge-Amp\`ere operator (see also \cite{Lu15}, \cite{ACC12} for domains). More precisely, we have.
 
\begin{lemma}   \label{Primitive} 
1) For any smooth path $t \longmapsto \phi_t$ in $\mathcal{E}^1 (\Omega)$ defined in some interval $I \subset \R$, we have
\begin{equation}
\frac{d}{d t} E (\phi_t) =  \int_\Omega (-\dot{\phi}_t) (dd^c \phi_t)^n, \, t \in I.
\end{equation}
In particular if $u, v \in \mathcal{E}^1 (\Omega)$ and $u \leq v$, then  $0 \leq E (v) \leq E (u)$. 

Moreover we have 
\begin{eqnarray*}
\frac{d^2}{d t} E (\phi_t) & = & \int_\Omega (- \ddot{\phi}_t) (dd^c \phi_t)^n   \\
& + &  n  \int_\Omega d \, \dot{\phi}_t \wedge d^c \dot{\phi}_t \wedge (dd^c \phi_t)^{n - 1}, \, t \in I
\end{eqnarray*}
In particular, for any  $u, v \in \mathcal{E}^1 (\Omega)$,  the function $t \longmapsto E( (1-t) u + t v)$ is a convex function in $[0,1]$.

2) The functional $E : \mathcal{E}^1 (\Omega) \longrightarrow \R^+$ is lower semi-continuous on $\mathcal{E}^1 (\Omega)$ for the $L^1_{loc}(\Omega)$-topology. It's continuous along decreasing sequences converging in $\mathcal{E}^1 (\Omega)$.

\end{lemma}
From this result it follows that if $u \in  \mathcal{E}^1 (\Omega)$ then
$$
E(u) = \inf \{E(\phi) \, ; \, \phi \in \mathcal{E}^0 (\Omega), \, \phi  \geq u\}\cdot
$$
It turns out that $ \mathcal{E}^1 (\Omega)$ is the natural domain of definition of the functional $E$.

\subsection{Integrability properties}

Let us  prove the second Claim used in the previous section (see \cite{BZ23}).
  \begin{lemma} \label{lem:integrability} Under the above assumptions, we have the following properties:
  
  $(i)$ $\mathcal E^1 (\Omega) \subset L^{n+1}(\Omega,\mu), $
  
  $(ii)$  $\phi \in \mathcal E^1 (\Omega) \setminus \{0\}  \Longrightarrow 0 < I_\mu(\phi) < +\infty.$
  \end{lemma}   
  \begin{proof}  Let $\rho$ be a negative strictly plurisbharmonic function in $\Omega$ such that  $dd^c \rho \geq \beta$ in $\Omega$. Since $f $ is bounded in $\Omega$, we can choose a positive constant $A \geq  \Vert f\Vert_{L^{\infty}(\Omega)}$ so that $(dd^c (A \rho))^n \geq f  d V = \mu $ in $\Omega$.
  Then by integration by parts, we obtain that for any $\phi \in \mathcal E^0(\Omega)$,
  $$
  \int_\Omega (-\phi)^{n+1} d \mu \leq A^n  \int_\Omega (-\phi)^{n+1} (dd^c \rho)^n \leq n ! A^n \Vert \rho\Vert_{L^{\infty}(\Omega)}^n \int_\Omega (-\phi) (dd^c \phi)^n.
  $$
  By the approximation Theorem of Cegrell (\cite{Ceg98}), the previous inequality holds for any $\phi \in \mathcal E^1 (\Omega)$, hence $\mathcal E^1 (\Omega) \subset L^{n+1}(\Omega,\mu)$.
  
  Now let $\phi \in \mathcal E^1 (\Omega) \setminus \{0\}$. Let $K \subset \Omega$ be a compact set such that  $ \int_K f d V > 0$. Then $M_K := max_K \phi < 0$ and 
  $$
  \int_\Omega (-\phi)^{n+1} d \mu \geq (- M_K) \int_K f d V > 0.
  $$
  \end{proof}
\subsection{Monotonicity property } 
Recall that for $\phi \in \mathcal E^1 (\Omega)\setminus \{0\}$, we have
 $$ 
 R(\phi) = \frac{E(\phi)}{\Vert \phi\Vert^{n+1}_{L^{n+1}(\Omega,\mu)}},
 $$  
 and then
$$
R(\phi) \Vert \phi\Vert_{L^{n+1}(\Omega,\mu)}^{n} = \frac{E(\phi)}{\Vert \phi\Vert_{L^{n+1}(\Omega,\mu)}}\cdot
$$

\begin{lemma}  \label{lem:decrease}
For any $k \in \N$, we have

\begin{equation}
\frac{E(u_{k+1})}{ \Vert u_{k+1} \Vert_{L^{n+1} (\Omega,\mu)}} \leq \frac{E(u_{k})}{ \Vert u_{k} \Vert_{L^{n+1} (\Omega,\mu)}}\cdot
\end{equation}
or equivalently
\begin{equation}
R(u_{k+1}) \Vert u_{k+1} \Vert^{n}_{L^{n+1} (\Omega,\mu)} \leq R(u_{k}) \Vert u_{k} \Vert^n_{L^{n+1} (\Omega,\mu)}.
\end{equation}
\end{lemma}
\begin{proof}  Multiplying the equation  $(dd^c u_{k+1})^n = R(u_k) (-u_k)^n \mu$ by $- u_{k+1}$ and integrating yields 
$$
E(u_{k+1}) = \int_\Omega (-u_{k+1}) (dd^c u_{k+1})^n = R(u_k) \int_\Omega (-u_{k+1}) (-u_k)^n d \mu.
$$

Applying H\"older inequality with the exponent $p=n +1$ and and its conjugate $q = (n+1)\slash n$, we obtain
$$
 \int_\Omega (-u_{k+1}) (-u_k)^n d \mu   \leq \Vert u_{k+1}\Vert_{L^{n+1}(\Omega,\mu)} \Vert u_{k}\Vert_{L^{n+1}(\Omega,\mu)}^n.
$$

Using the definition of $R(u_k)$ we obtain the required inequality
$$
E(u_{k+1})  \leq E(u_k) \frac{\Vert u_{k+1}\Vert_{L^{n+1}(\Omega,\mu)}}{ \Vert u_{k}\Vert_{L^{n +1}(\Omega,\mu)}}\cdot
$$
\end{proof}
\subsection{Uniform bounds on the Rayleigh quotients}
Let $(\lambda_1,w) $ be the solution to \eqref{eq:VP-MA} normalized by the condition  $\Vert w \Vert_{L^{\infty}(\Omega)} = \lambda_1^{-1}$ given by Theorem A.

\begin{corollary} \label{eq:Rayleigh bound} For any $k \in \N$, we have 

$(i)$ $u_k \leq w < 0$ in $\Omega$, hence $0 < \Vert w \Vert_{L^{n+1} (\Omega,\mu)} \leq \Vert u_k\Vert_{L^{n+1} (\Omega,\mu)}$.
\smallskip

$(ii)$ $R(u_k) \leq C_0 := R(u_{0}) \Vert u_{0} \Vert^n_{L^{n+1}(\Omega,\mu)}  (\Vert w \Vert_{L^{n+1} (\Omega,\mu)})^{- n}$.
\end{corollary}

\begin{proof} We prove by induction that $u_k \leq w$ in $\Omega$ for any $k \in \N$.

Indeed since $ 0 \leq -\lambda_1 w \leq 1$ in $\Omega$, it follows that 
$$
(dd^c w)^n  = (-\lambda_1 w)^n \mu \leq \mu\leq (dd^c u_0)^n.
$$

Since $u_0 \leq 0 = w$ in $\partial \Omega$, by the Comparison Principle \cite{BT76},  $u_0 \leq w$ in $\Omega$.

Now assume that $u_k \leq w$ for some $k \in \N^*$,  then 
\begin{eqnarray*}
(dd^c w)^n & = & (-\lambda_1 w)^n \mu \\
& \leq & \lambda_1^n (-u_k)^n \mu
 \leq R(u_k) (-u_k)^n \mu  
= (dd^c u_{k+1})^n.
\end{eqnarray*}

Since $u_{k+1} =0= w $ in $\partial \Omega$, it follows from the Comparison Principle  that $u_{k+1} \leq w < 0$ in $\Omega$. 
This inequality implies that for any $k \in \N$,
$$\Vert u_{k} \Vert_{L^{n+1} (\Omega,\mu)}^n \geq \Vert w \Vert_{L^{n+1} (\Omega,\mu)}^n  > 0.$$

\smallskip

By  the monotonicity property Lemma \ref{lem:decrease}, we deduce that for any $k \in \N$,
$$
R(u_{k}) \leq C_0 := R(u_{0}) \Vert u_{0} \Vert^n_{L^{n+1}(\Omega,\mu)}  (\Vert w \Vert_{L^{n+1} (\Omega,\mu)})^{- n}.
$$
\end{proof}

\subsection{ Uniform H\"older estimates}
Recall that for each $k\in \N$, 
$$
(dd^c u_{k +1})^n = f_k d V,  \, \, \text{where} \, \,  f_k := R(u_k) (-u_k)^n f. 
$$ 
 It follows from the monotonicity property that for any $k \in \N^*$,
\begin{eqnarray*}
\Vert f_k\Vert_{L^{1 + 1 \slash n}(\Omega)} &\leq& \Vert f\Vert_{L^{\infty}(\Omega)}^{1 \slash (n+1)} R(u_k)  \left(\int_\Omega (-u_k)^{n+1} f d V\right)^{n \slash (n+1)} \\ 
& \leq & \Vert f\Vert_{L^{\infty}(\Omega)}^{1 \slash (n+1)}  R(u_k) \Vert u_k \Vert_{L^{n+1} (\Omega,\mu)}^n\\ 
& \leq &\Vert f\Vert_{L^{\infty}(\Omega)}^{1 \slash (n+1)} R(u_0) \Vert u_0 \Vert_{L^{n+1} (\Omega,\mu)}^n =: C'_0
\end{eqnarray*}
 H\"older  uniform a priori estimates imply : $\exists \, \,   \alpha \in ]0,1[$,  $\exists \,  M = M(n,\alpha,\Omega,C'_0)  > 0$,  such that $ \forall k \in \N :$
\begin{equation} \label{eq:HolderEst}  
\Vert u_{k+1} \Vert_{L^{\infty}(\Omega)}  \leq \Vert u_{k+1} \Vert_{C^{\alpha}(\Omega)} \leq M.
\end{equation}

The uniform estimate is proved in \cite{Kol96}, while the H\"older estimate is proved in  \cite{GKZ08} in a particular case and  in  \cite{Ch15} in the general case. Actually it's show that we can take any $0 < \alpha < \frac{2}{n (n+1) + 1}\cdot$ 
  
\section{Proof of the Main Theorem}

The proof will be done in several steps following the same scheme as \cite{AK20}.
\subsection{Existence of a limit point}  Thanks to the uniform estimates \eqref{eq:HolderEst}, we can apply Arzel\`a-Ascoli theorem to obtain an increasing sequence $(k(j))_{j \in \N}$ in $\N^*$  such that
 \begin{itemize} 
 \item $(u_{k(j)})_{j\in \N} \to \varphi \in PSH(\Omega) \cap C^{\alpha}(\bar \Omega)$  uniformly on $\bar \Omega$ as $j \to +\infty$, and 
 \item $(u_{k(j)+1})_{j \in \N} \to \psi \in PSH(\Omega) \cap C^{\alpha}(\bar \Omega)$  uniformly on $\bar \Omega$ as $j \to +\infty$.  
 \end{itemize}
 \smallskip

The goal is to prove that $\varphi = \psi$.
  We will use the  following elementary fact. 
 
 \smallskip
 \begin{lemma} \label{lem:Energyconvergence} Let $(\phi_j)$ be a sequence in $PSH_0(\Omega)  \cap C^0(\bar \Omega)$ which converges uniformly in $\bar \Omega$ to a function $\phi \in PSH_0 (\Omega)  \cap C^0(\bar \Omega)$. Assume that there exists a bounded Borel measure $\nu$  on $\Omega$ such that $\mathcal E^1(\Omega) \subset L^1(\Omega,\nu)$ and  for any $j \in \N$, $(dd^c \phi_j)^n \leq \nu$ on $\Omega$. Then we have
\begin{itemize}
\item $\lim_{j \to + \infty} E (\phi_j) = E(\phi)$;
\item  $ \lim_{j \to + \infty} I_\mu (\phi_j) = I_\mu (\phi)$;  
\item  $ \lim_{j \to + \infty} R (\phi_j) = R(\phi)$.
\end{itemize}
\end{lemma}
\begin{proof} It's enough to prove that $\lim_{j\to + \infty} E(\phi_j) = E(\phi)$.
 We have for $j \in \N$, 
\begin{eqnarray*}
E (\phi_j) - E (\phi) &=& \int_\Omega (\phi - \phi_j) (dd^c \phi_j)^n + \int_\Omega (- \phi) \left[(dd^c \phi_j)^n - (dd^c \phi)^n\right]\\
&=& A_j + B_j.
\end{eqnarray*}
It follows from our assumptions that for any $j \in \N$, we have
\begin{eqnarray*}
\vert A_j \vert = \left\vert \int_\Omega (\phi - \phi_j) (dd^c \phi_j)^n \right\vert  &\leq & \Vert\phi_j - \phi\Vert_{C^0(\bar \Omega)}  \nu (\Omega), 
\end{eqnarray*}
which proves that $\lim_{j \to + \infty} A_j = 0$.

On the other hand, fix a compact set $K \subset \Omega$.  Then for any $j \in \N$, we have
\begin{eqnarray} \label{eq:decomposition}
 B_j  &=&  \int_\Omega (- \phi) \left[(dd^c \phi_j)^n - (dd^c \phi)^n\right]  \nonumber\\
 &=&  \int_{\Omega\setminus K} (-\phi)\left[(dd^c \phi_j)^n - (dd^c \phi)^n\right] +   \int_K (-\phi) \left[(dd^c \phi_j)^n - (dd^c \phi)^n\right].
\end{eqnarray}
The first term of the RHS of the formula \eqref{eq:decomposition} is dominated by $2 \int_{\Omega\setminus K} (- \phi) d \nu$, which can be made arbitrarely small by choosing the compact set $K$ large enough, since $\int_\Omega (-\phi) d \nu < + \infty$.

We claim that the second term converges to $0$ as $j \to+ \infty$. Indeed,
observe that $(-\phi)(dd^c \phi_j)^n \to (-\phi) (dd^c \phi)^n$ in the sense of currents on $\Omega$. Then by choosing the compact set $K$ so that $\int_{\partial K} (-\phi)  (dd^c \phi)^n = 0$, we deduce that 
$$
\lim_{j \to + \infty}  \int_K (-\phi) (dd^c \phi_j)^n =  \int_K (-\phi) (dd^c \phi)^n.
$$
Therefore the second term of the RHS of the formula \eqref{eq:decomposition} converges to $0$ as $j \to + \infty$.
This proves that $\lim_{j\to + \infty} B_j = 0$ and then $\lim_{j\to + \infty} E(\phi_j) = E(\phi)$, which proves the lemma.
\end{proof}

Next we prove the following result. 
\begin{lemma} \label{lem:Conv1} We have  $\psi = \varphi$ in $\Omega$ and $R(\varphi) = \lambda_1^n$. 
\end{lemma}
\begin{proof}
By the uniform estimates \eqref{eq:Rayleigh bound} and \eqref{eq:HolderEst}, we see  that there exists a uniform constant $C > 0$ such that for any $j\in \N$, we have
$$
(dd^c u_{k(j)})^n = R(u_{k(j)-1}) (-u_{k(j)-1})^n f  d V \leq C d V.
$$ 

Applying Lemma \ref{lem:Energyconvergence}, we conclude that
 \begin{itemize}
 \item $\lim_{j \to + \infty} E (u_{k(j)}) = E(\varphi)$ and $\lim_{j \to + \infty} I_\mu (u_{k(j)}) = I_\mu (\varphi)$
 \item $\lim_{j \to + \infty} R (u_{k(j)}) = R (\varphi)$.
  \end{itemize}
  
  \smallskip
 
Recall that for any $ j \in \N$,  $(dd^c u_{k(j) +1})^n = R (u_{k(j)}) (-u_{k(j)})^n f d V$. Taking the limit  we obtain { $(dd^c \psi)^n= R(\varphi)(-\varphi)^n f d V$} on $\Omega$.

\smallskip

Now  since $k(j+1) \geq k(j)+1 \geq k(j)$, the monotonicity property Lemma \ref{lem:decrease} implies that for any $j \in \N$,
\begin{eqnarray*}
R(u_{k(j+1)}) \Vert u_{k(j+1)}\Vert_{L^{n+1}(\Omega,\mu_f)}^n  &\leq & R(u_{k(j)+1}) \Vert u_{k(j)+1}\Vert_{L^{n+1}(\Omega,\mu)}^n \\
& \leq & R(u_{k(j)}) \Vert u_{k(j)}\Vert_{L^{n+1}(\Omega,\mu)}^n.
\end{eqnarray*}

Taking the limit, we obtain {$R(\psi) \Vert \psi\Vert_{L^{n+1}(\Omega,\mu)}^n = R(\varphi) \Vert \varphi\Vert_{L^{n+1}(\Omega,\mu)}^n.$}

On the other hand, since $(dd^c \psi)^n = R(\varphi) (-\varphi)^n d \mu$, applying H\"older inequality, we obtain
\begin{eqnarray*}
R(\psi) \Vert \psi\Vert_{L^{n+1}(\Omega,\mu)}^{n+1} &= & E(\psi) = \int_\Omega (-\psi) (dd^c \psi)^n 
=  R(\varphi) \int_\Omega (-\psi) (-\varphi)^n d \mu \\ 
&\leq &  R(\varphi) \Vert \psi\Vert_{L^{n+1}(\Omega,\mu)}  \Vert \varphi\Vert_{L^{n+1}(\Omega,\mu)}^n
= R(\psi) \Vert \psi\Vert_{L^{n+1}(\Omega,\mu)}^{n+1}.
\end{eqnarray*}

This means that we have equality in the H\"older inequality, hence there exists a constant $c > 0$ such that
$(-\psi)^{n +1}= c (-\varphi)^{n+1}$ in $\Omega$ i.e. $\psi = c^{1\slash (n+1)} \varphi$. 

Previous inequalities imply $c= 1$, hence $\psi = \varphi$.

\smallskip
 Summarizing we see that  {$(u_{k(j)})_{j \in \N} \longrightarrow \varphi \in PSH(\Omega) \cap C^\alpha(\bar \Omega)$} uniformly in $\bar \Omega$ and {$(dd^c \varphi)^n = R(\varphi) (-\varphi)^n f dV$} with $\varphi = 0$ in $\partial \Omega$ and $\varphi \leq w<0$.

Hence {$R(\varphi) = \lambda_1^n$} by uniqueness of the eigenvalue in Theorem A.
\end{proof}

 \subsection{Convergence of the sequence  $(u_k)$} 
  Here we prove the uniqueness of the limit point.
  \begin{lemma} \label{lem:Conv2} The sequence $(u_k)$ converges uniformly on $\bar \Omega$ to $\varphi$.
  
  \end{lemma}
  \begin{proof}
 Assume there are two subsequences $(u_{k(j)}) \to \varphi$  and  $(u_{k' (j)}) \to \varphi'$   uniformly on $\bar \Omega$.

 By  Lemma \ref{lem:Conv1}, we have $R(\varphi) =  \lambda_1^n=R(\varphi')$. By the uniqueness property of Theorem A we conclude that {$\varphi = c \cdot \varphi'$} in $\Omega$ for some constant $c > 0$.

 It remains to prove that $c = 1$.

 Indeed mixing the two sequences, it is possible to construct a subsequence $(\ell(j))_{j \in \N}$ of $(k(j))_{j \in \N}$ and a subsequence  $(\ell'(j))_{j \in \N}$ of $(k'(j))_{j \in \N}$ such that $\ell(j) < \ell'(j) < \ell(j+1)$ for any $j\in \N$.
 
 \smallskip
 
  Then by the monotonicity property, we have
  \begin{eqnarray*}
 R(u_{\ell(j+1)}) \Vert u_{\ell(j+1)}\Vert_{L^{n+1}(\Omega,\mu)}^{n} &\leq & R(u_{\ell'(j)}) \Vert u_{\ell'(j)}\Vert_{L^{n+1}(\Omega,\mu)}^{n} \\
 & \leq & R(u_{\ell(j)}) \Vert u_{\ell(j)}\Vert_{L^{n+1}(\Omega,\mu)}^{n}\cdot
  \end{eqnarray*}
  
  Taking the limits we obtain $R(\varphi) \Vert \varphi\Vert_{L^{n+1}(\Omega,\mu)}^{n} = R(\varphi') \Vert \varphi'\Vert_{L^{n+1}(\Omega,\mu)}^{n}$, hence $\Vert \varphi\Vert_{L^{n+1}(\Omega,\mu)}^{n} = \Vert \varphi'\Vert_{L^{n+1}(\Omega,\mu)}^{n}$.
  The uniqueness property  of a weak solution provided by Theorem A yields $\varphi = \varphi'$.
The  Main Theorem follows then from Lemma \ref{lem:Conv1} and Lemma \ref{lem:Conv2}.
\end{proof}

\smallskip

To conclude this note we address the following open question.
\smallskip
\smallskip

\noindent{\bf Question :}
 {\it Assume that we choose the initial function $u_0$ so that $(dd^c u_0)^n = f d V$ on $\Omega$ and ${u_0}_{\mid \partial \Omega} = 0$.  We have seen that the sequence $(u_k)_{k \in \N}$ converges uniformly to an eigenfunction $\varphi \in PSH(\Omega) \cap C^0(\bar \Omega)$ such that $\varphi \leq w$ in $\Omega$, where  $ w $ is the eigenfunction  normalized by the condition $\Vert w \Vert_{L^{\infty}(\Omega)} = \lambda_1^{-1}$. 
 
 Is it possible  to characterize the limit function $\varphi$ ? 
 Is it true that $\varphi = w$?}
 
 \smallskip
 \smallskip
 \smallskip
 
 \noindent{\bf Aknowledgements :}  It's a great pleasure for me to dedicate this article to my friend Wies\l aw Plesniak on the occasion of this $80^{\text{th}}$ birthday. 
 
  Wies\l aw  has made important contributions  to  Pluripotential Theory and the Theory of Approximation who inspired many mathematicians. In particular I am personally grateful to him for introducing me to the fascinating world of Markov Inequalities and their applications.  
 
 Finally the author would like to thank the anonymous reviewers for their careful reading of the first version of this article, which helped to correct many typos.


\begin{thebibliography}{widestlabel}

\bibitem[AK20]{AK20}  F. Abedin, J. Kitagawa : {\em Inverse iteration for the Monge-Amp\`ere eigenvalue problem}. Proceedings AMS, 148 (2020), no 11, 4875-4886. 

\bibitem[ACC12]{ACC12} P. Åhag,  U. Cegrell, R.  Czyz :  {\it On Dirichlet's principle and problem.} Math. Scand. 110 (2), 235-250 (2012).
	
  \bibitem[BZ23]{BZ23} P. Badiane, A. Zeriahi : {\it  The eigenvalue problem for the complex Monge-Amp\`ere opertaor}. J. Geom. An. (2023), 33:367, 44pp.
  
   \bibitem[BT76]{BT76} E. Bedford, B. A. Taylor : {\it The Dirichlet problem for a complex Monge-Amp\`ere equation.} Invent. Math., 37 (1), 1-44 (1976).	
	
 \bibitem[BT82]{BT82} E. Bedford, B. A. Taylor :	{\it A new capacity for plurisubharmonic functions.} Acta Math. 149 , no. 1 (2), 1-40 (1982).
 

  \bibitem[BNV94]{BNV94} H. Berestycki, L. Nirenberg, S.R.S. Varadhan : {\it   The first eigenvalue and maximum principle for second order elliptic differential operators in general domains.} Comm. Pure Appl. Math. 47 (1), 47-92 (1994).


 \bibitem[BBGZ13]{BBGZ13} R. Berman, S. Boucksom, V. Guedj, A. Zeriahi : {\it A variational approach to complex Monge-Amp\`ere equations.} Publ. Math. Inst. Hautes Études Sci. 117, 179-245  (2013).
	
 \bibitem[CKNS85]{CKNS85} L. Caffarelli, J. J. Kohn, L. Nirenberg, J. Spruck : {\it The Dirichlet problem for nonlinear second-order elliptic equations. II. Complex Monge-Amp\`ere, and uniformly elliptic, equations.} Comm. Pure Appl. Math. 38 (2), 209-252 (1985).

 \bibitem[Ceg98]{Ceg98}  U. Cegrell : {\it Pluricomplex energy.} Acta Math. 180 (2), 187-217  (1998). 
	
	
 \bibitem[Ch15]{Ch15} M. Charabati : {\it Hölder regularity for solutions to complex Monge-Amp\`ere equations.}
Ann. Polon. Math. 113 (2), 109-127  (2015).

\bibitem[CLMcC24]{CLMcC24} J. Chu, Y. Liu,  N. McCleerey : {\it The eigenvalue problem for the complex Hessian opertaor on $m$psuedoconvex manifolds}. Preprint arXiv:2402.03098v1.


 \bibitem[Eva10]{Eva10} L. C. Evans : {\it Partial Differential Equations.} Graduate Studies in Mathematics 19, Second edition, American Mathematical Society (2010).
	
\bibitem[Gav77]{Gav77} B. Gaveau : {\it  M\'ethode de contr\^ole optimal en analyse complexe, I.} J. Funct. Anal. 25,  391-411 (1977).
	
 \bibitem[GKZ08]{GKZ08} V. Guedj, S. Koldziej, A. Zeriahi : {\it Hölder continuous solutions to Monge-Amp\`ere equations.} Bull. Lond. Math. Soc. 40 (6), 1070–1080  (2008). 
	
 \bibitem[GZ17]{GZ17} V. Guedj, A. Zeriahi : {\it Degenerate Complexe Monge-Amp\`ere Equations.} EMS Tracts in Mathematics 26 (2017).
 
\bibitem[Kol96]{Kol96}  S. Kolodziej : {\it  Some sufficient conditions for solvability of the Dirichlet problem for the complex Monge-Amp\`ere operator}.  Ann. Pol. Math. 65 (1), p. 11-21  (1996).
	
 \bibitem[Kol03]{Kol03}  S. Kolodziej : {\it Equicontinuity of families of plurisubharmonic functions with bounds on their Monge-Amp\`ere masses.} Math. Z. 240 (4), 835-847  (2002).
	
 
 \bibitem[Lions86]{Lions86} P. L. Lions : {\it Two remarks on the Monge-Amp\`ere equations.} Ann. Mat. Pura Appl. 142 (4), 263-275 (1986).
	
 \bibitem[Lu15]{Lu15} C. H. Lu : {\it  A variational approach to complex Hessian equations in $\C^n$.} J. Math. Anal. Appl. 431 (1), 228-259  (2015).	
	
	
 \bibitem[Tso90]{Tso90} K. Tso : {\it On a real Monge-Amp\`ere functional.} Invent. Math. 101, 425-448 (1990).
	
\end{thebibliography}
\end{document}